\documentclass[10pt]{article}

\usepackage{amsmath}
\usepackage{amssymb}  
\usepackage[T1]{fontenc}
\usepackage{pifont}
\usepackage{pstricks}
\usepackage{amsfonts}
\usepackage{graphicx}
\usepackage{amsmath}
\usepackage{authblk}
\usepackage{pstricks} 
\usepackage{pstricks-add}
\usepackage{caption}
\usepackage{placeins}
\usepackage{pstricks}
\usepackage{pdfpages}
\usepackage{framed}
\usepackage{cancel}
\usepackage{bm}
\usepackage{epstopdf}
\def\x2dot{\mathop{x}\limits}
\def\y2dot{\mathop{y}\limits}

\def\bfy2dot{\mathop{\bf y}\limits}
\def\z2dot{\mathop{z}\limits}

\def\csi2dot{\mathop{\xi}\limits}
\def\et2dot{\mathop{\eta}\limits}
\def\bet2dot{\mathop{\beta}\limits}
\def\t2dot{\mathop{\theta}\limits}
\def\s2dot{\mathop{\sigma}\limits}
\def\d2dot{\mathop{\delta}\limits}
\def\q2dot{\mathop{q}\limits}
\def\l2dot{\mathop{\lambda}\limits}
\def\ps2dot{\mathop{{\cal E}}\limits}
\def\tet2dot{\mathop{\theta}\limits}
\def\bfx2dot{\mathop{\bf X}\limits}
\def\bfy2dot{\mathop{\bf y}\limits}
\def\bfq2dot{\mathop{\bf q}\limits}
\def\bbfq2dot{\mathop{\bar {\bf q}}\limits}
\def\w2{\mathop{W}\limits}
\def\xgrande2dot{\mathop{\bf X}\limits}
\def\p02dot{\mathop{P}\limits}
\def\a2dot{\mathop{A}\limits}

\newtheorem{lem}{Lemma}
\newtheorem{prop}{Proposition}
\newtheorem{cor}{Corollary}

\newtheorem{rem}{Remark}
\newtheorem{exe}{Example}

\title{Energy balance for nonholonomic nonlinear systems}
\author{F.~Talamucci}
\affil{{\it DIMAI, Dipartimento di Matematica e Informatica ``Ulisse Dini''},\\
{\it	Universit\`a degli Studi di Firenze, Italy}\\
{\it	e-mail: federico.talamucci@unifi.it}}
\date{}

\begin{document}
	\bibliographystyle{plain}
	
	\setcounter{equation}{0}

	\maketitle
	
	\vspace{.5truecm}
	
	\noindent
	{\bf 2010 Mathematics Subject Classification:} 37J60, 70F25, 70H03.
	
	\vspace{.5truecm}
	
	\noindent
	{\bf Keywords:} Lagrange's equations of motion - Nonholonomic systems - Nonlinear kinematic constraints - Equations of Voronec type - Generalized Theorem of energy.
	
	\vspace{.5truecm}
	
	\noindent
	{\bf Abstract}. We consider nonholonomic systems with nonlinear restrictions with respect to the velocities. The mathematical problem is formulated by means of the Voronec equations extended to the nonlinear case. The main point of the paper is the balance of the mechanical energy induced by the equations of motion; the conservation of the energy on the basis of the tipology of the constraint equations is discussed. Several examples are performed.

\section{Introduction}

\noindent
Nonholonomic systems are frequently encountered in mechanical and engineering problems and attention 
to mathematical models concerning kinematic constraints 
is increasingly paid. 

\noindent
The  trait ``nonholonomic'' refers, in a general way, to restrictions imposed on a system which are not 
expressible solely in terms of one or more equations involving only the spatial coordinates which delineate the position of the system.

\noindent
In the wide range of restrictions with nonholonomic features, we focus the interest on the 
kinematic constraints, that is restrictions which can be formulated by equations containing the coordinates and components of the velocities. As a matter of fact, this kind of constraints show 
a more suitable relevance and evidence concerning applications and feasibility.

\noindent
A remarkable difference exists between the case of linear kinematic constraints and the case of nonlinear ones: even in literature, the starting point of the theory of linear constraints is distant in time (\cite{capligin}, \cite{hamel}) and in parallel to the development of the lagrangian formalism, 
whereas the advancements in nonlinear constrained systems are somewhat recent.
The not uniform progress is motivated both by modellistic reasons (actually, the linear case can be  convincingly considered as a natural extent of the holonomic case, as we will see just below)
and by practical reasons of realization of physical models exhibiting nonlinear constraints (such a question, arised in \cite{appell} and studied in \cite{zek1}, is still debated, see \cite{benenti}).

\noindent
In this paper we deal with systems submitted to nonlinear nonholonomic constraints.
The first part (Section 2) consists in recalling the equations of motions proposed in a previous work (\cite{tal}) and based on an elementary principle, formulated along the lines of the linear nonholonomic case. The main requirements of this formulations are the explicit form of the constraint equations and the
use of some of the generalized velocities as independent parameters. In this way, the formal path is simple and proceeds with extending the classical theory of mechanics.
We will also point out the correspondence of the obtained equations of motion with models on the same topics present in literature, where the equations derive from more refined and formally complex methods.
In literature, nonholonomic mechanics appear to be studied mainly via geometric methods, as Lagrangian systems on fibered manifolds (see, among others,  \cite{leon}): the method proposed in the present paper is a handly approach developing in analogy to the basic concepts of a holonomic system.

\noindent
The second part (Section 3) is devoted to the main result, regarding the expression for the rate of change of the mechanical energy of the system. 
The selected typology of the equations of motion turns ut to be appropriate both from the mathematical point of view and for the physical intepretation, since the variables involved are the real velocities.

\noindent
The energy balance is achieved via a somehow usual handling of the equations of motion: in spite of the spontaneity of the method, the resulting formula shows distintly each contribute in terms of energies and it is suitable in order to investigate which categories of constraints 
entail the conservation of generalized energy (Jacobi's first integral). In Section 4 we identify some cases of nonholonomic systems where the energy balance equation actually infers the conservation of the mechanical energy.

\noindent
Finally, various examples aimed at the implementation of the results 
on both linear and nonlinear nonholonomic constrained systems have been added, by drawing basically the most recurrent models in literature.

\section{Mathematical model and notations}

\noindent
The starting point is a holonomic system ${\bf X}(q_1,\dots, q_n,t)$, where ${\bf X}$ locates the position of $N$ material points $(P_i,M_i)$, $i=1,\dots, N$. In addition to the $3N-n$ holonomic conditions, the system undergoes $k<n$ kinematic constraints

\begin{equation}
\label{constr}
\Phi_\nu(q_1,\dots, q_n, {\dot q}_1, \dots, {\dot q}_n, t)=0, \quad \nu=1,\dots, k
\end{equation}
which can be linear or nonlinear and are assumed to be independent, in the sense that the rank of the Jacobian matrix 
$J_{({\dot q}_1, \dots, {\dot q}_n)}(\Phi_1, \dots, \Phi_k)$ is the greatest value $k$.

\noindent
Under that assumption, if for istance $det\,J_{({\dot q}_{m+1}, \dots, {\dot q}_n)}(\Phi_1, \dots, \Phi_k)\not =0$, $m=n-k$, conditions (\ref{constr}) can be made explicit:
\begin{equation}
\label{constrexpl}
{\dot q}_{m+\nu}=\alpha_\nu(q_1, \dots, q_n, {\dot q}_1, \dots, {\dot q}_m, t), \quad \nu=1,\dots, k.
\end{equation}
with $m=n-k$. The case when all of (\ref{constr}) are linear is denoted by the functions $\alpha_{\nu,j}$ such that
\begin{equation}
\label{constrexpllin}
\alpha_\nu=\sum\limits_{j=1}^m a_{\nu,j}(q_1,\dots, q_n,t){\dot q}_j, \quad \nu=1,\dots, k.
\end{equation}

\subsection{Equations of motion}

\noindent
Let ${\bf X}^{(\mathsf{M})}\in {\Bbb R}^{3N}$ list the $n$ triplets $M_iP_i$, $i=1,\dots, N$, so that 
${\cal Q}={\dot {\bf X}}^{(\mathsf{M})}$ is the $3N$--vector of linear momenta; furthermore, let us denote by  ${\mathcal F}=({\mathcal F}_1, \dots, {\mathcal F}_N)$, ${\mathcal R}=({\mathcal R}_1,\dots,{\mathcal R}_N)$ the list in ${\Bbb R}^{3N}$ of the active forces ${\mathcal F}_i$ and of the constraint forces ${\mathcal R}_i$ concerning each $P_i$, $i=1,\dots, N$.

\noindent
We present initially the equations of motion in the newtonian form. The key statement is 
\begin{equation}
\label{ge}
\left( {\dot {\mathcal Q}}-{\mathcal F}- {\mathcal R}\right)\cdot {\widehat {\dot {\bf X}}}=0
\end{equation}
where 
\begin{equation}
\label{hatx}
{\widehat {\dot {\bf X}}}=
\sum\limits_{r=1}^m  {\dot q}_r\left(
\dfrac{\partial {\bf X}}{\partial q_r}+
\sum\limits_{\nu=1}^{k}\dfrac{\partial \alpha_\nu}{\partial {\dot q}_r}
\dfrac{\partial {\bf X}}{\partial q_{m+\nu}}\right)
\end{equation}
are all the possible displacements (see \cite{benenti}, \cite{tal}) consistent with the kinematic restrictions 
(\ref{constrexpl}),  for arbitrary  $({\dot q}_1, \dots, {\dot q}_m)\in {\Bbb R}^m$.  

\noindent
At this point, the assumptions of ideal constraints demands ${\cal R}\cdot {\hat {\dot {\bf X}}}=0$ for all possible (\ref{hatx}), so that (\ref{ge}) is equivalent to the equations
\begin{equation}
\label{ne}
\left( {\dot {\mathcal Q}}-{\mathcal F}\right)\cdot {\bf X}_r=0, \quad r=1, \dots, m, \quad {\bf X}_r=
\dfrac{\partial {\bf X}}{\partial q_r}+
\sum\limits_{\nu=1}^{k}\dfrac{\partial \alpha_\nu}{\partial {\dot q}_r}
\dfrac{\partial {\bf X}}{\partial q_{m+\nu}}
\end{equation}
which correspond explicitly to the $m$ differential equations (which have to be joined with (\ref{constrexpl}))

\begin{equation}
\label{vnl2}
\sum\limits_{r=1}^m \left( C_i^r {\q2dot^{..}}_r+
\sum\limits_{s=1}^m D_i^{r,s} {\dot q}_r{\dot q}_s + E_i^r {\dot q}_r\right) +G_i={\cal F}^{(q_i)}+
\sum\limits_{\nu=1}^k \dfrac{\partial \alpha_\nu}{\partial {\dot q}_i} {\cal F}^{(q_{m+\nu})}, 
\qquad i=1,\dots, m
\end{equation}
where $\mathcal{F}^{(q_j)}=\mathcal{F}\cdot \dfrac{\partial {\bf X}}{\partial q_j}$ for each $j=1, \dots, n$ is the generalized force, 
$C_i^r$, $D_i^{r,s}$, $E_i^r$ and $G_i$, $i,r,s=1,\dots, m$, depend on $q_1, \dots, q_n$, ${\dot q}_1$, $\dots$, ${\dot q}_m$, $t$ and are defined by

\begin{eqnarray}
\nonumber
C_i^r&=&g_{i,r}+\sum\limits_{\nu,\mu=1}^k\left(g_{i,m+\nu}\dfrac{\partial \alpha_\nu}{\partial {\dot q}_r}
+g_{m+\nu,r}\dfrac{\partial \alpha_\nu}{\partial {\dot q}_i}+ g_{m+\nu,m+\mu}\dfrac{\partial \alpha_\mu}{\partial {\dot q}_i}\dfrac{\partial \alpha_\nu}{\partial {\dot q}_r}\right), \quad
D_i^{r,s}=\xi_{r,s,i}+\sum\limits_{\nu=1}^k\xi_{r,s,m+\nu}\dfrac{\partial \alpha_\nu}{\partial {\dot q}_i}, \\
\label{coeff}
\nonumber
E_i^r &=&  \sum\limits_{\nu,\mu=1}^k \left(
2\xi_{r,m+\nu,m+\mu} \alpha_\nu + g_{m+\nu,m+\mu}\dfrac{\partial \alpha_\nu}{\partial q_r}+2 \eta_{r,m+\mu} \right)\dfrac{\partial \alpha_\mu}{\partial {\dot q}_i}+
\label{coeffcdef}
\sum\limits_{\nu=1}^k \left( 2\xi_{r,m+\nu,i} \alpha_\nu + 
g_{m+\nu, i}\dfrac{\partial \alpha_\nu}{\partial q_r} \right)
+2 \eta_{r,i}, \\
\nonumber
G_i&=&
\sum\limits_{\nu,\mu,p=1}^k
\left(  \xi_{m+\nu, m+\mu,m+p}\alpha_\nu\alpha_\mu + g_{m+\nu,m+p} \alpha_\mu 
\dfrac{\partial \alpha_\nu}{\partial q_{m+\mu}}+
2\eta_{m+\nu, m+p}\alpha_\nu +g_{m+\nu, m+p} \dfrac{\partial \alpha_\nu}{\partial t}+\zeta_{m+p}\right) \dfrac{\partial \alpha_p}{\partial {\dot q}_i} \\
\nonumber
&+&
\sum\limits_{\nu,\mu=1}^k\left(
\xi_{m+\nu,m+\mu,i}\alpha_\nu\alpha_\mu+g_{m+\nu,i}\alpha_\mu \dfrac{\partial \alpha_\nu}{\partial q_{m+\mu}}
\right) +
\sum\limits_{\nu=1}^k \left( 
2 \eta_{m+\nu,i} \alpha_\nu +g_{m+\nu,i} \dfrac{\partial \alpha_\nu}{\partial t}\right)+\zeta_i.
\end{eqnarray}
by setting, for any $i,j,k=1, \dots, n$:
\begin{equation}
\label{gxi}
\begin{array}{ll}
g_{i,j}(q_1,\dots, q_n, t)=\dfrac{\partial {\bf X}^{(\mathsf{M})}}{\partial q_i}\cdot 
\dfrac{\partial {\bf X}}{\partial q_j}, &
\xi_{i,j,k}(q_1, \dots, q_n, t)=\dfrac{\partial^2{\bf X}^{(\mathsf{M})}}{\partial q_i\partial q_j}\cdot \dfrac{\partial {\bf X}}{\partial q_k},\\
\\
\eta_{i,j}(q_1,\dots, q_n, t)= \dfrac{\partial^2 {\bf X}^{(\mathsf{M})}}{\partial q_i \partial t}\cdot \dfrac{\partial {\bf X}}{\partial q_j}, &
\zeta_i(q_1,\dots, q_n, t)=\dfrac{\partial^2 {\bf X}^{(\mathsf{M})}}{\partial t^2}\cdot 
\dfrac{\partial {\bf X}}{\partial q_j}.
\end{array} 
\end{equation}	

\noindent
The explicit form (\ref{ne}) is convenient in order to easily disclose the terms with the second derivatives 
${\q2dot^{..}}_r$, $r=1, \dots, m$, which appear only in the linear terms with coefficients $C_i^r$.
In case of scleronomic holonomous constraints ${\bf X}({\bf q})$, in (\ref{gxi}) it is $\eta_{i,j}=0$, $\zeta_i=0$
for any $i,j=1, \dots, n$. 

\begin{prop}
The $m\times m$ matrix $C$ with entries $C_i^r$, $r=1, \dots, m$, is positive definite.

\end{prop}
{\bf Proof}: 
Indeed the terms originate from
$C_i^r=
\left(\dfrac{\partial {\bf X}^{(M)}}{\partial q_r}+
\sum\limits_{\nu=1}^{k}\dfrac{\partial \alpha_\nu}{\partial {\dot q}_r}
\dfrac{\partial {\bf X}^{(M)}}{\partial q_{m+\nu}}\right)\cdot 
\left(
\dfrac{\partial {\bf X}}{\partial q_i}+
\sum\limits_{j=1}^{k}\dfrac{\partial \alpha_j}{\partial {\dot q}_i}
\dfrac{\partial {\bf X}}{\partial q_{m+j}}\right)$. 
Setting ${\bf X}_r^{(\sqrt{M})}=\dfrac{\partial {\bf X}^{(\sqrt{M})}}{\partial q_r}+
\sum\limits_{\nu=1}^{k}\dfrac{\partial \alpha_\nu}{\partial {\dot q}_r}
\dfrac{\partial {\bf X}^{(\sqrt{M})}}{\partial q_{m+\nu}}$, where ${\bf X}^{(\sqrt{M})}$ stands for $(\sqrt{m_1}P_1, \dots, \sqrt{m_n}P_n)$, one has $C_\nu^i= {\bf X}_\nu^{(\sqrt{M})}\cdot {\bf X}_i^{(\sqrt{M})}$, so that the matrix $C$ is positive definite, since the vectors ${\bf X}_j^{(\sqrt{M})}$, $j=1, \dots, m$ are linearly independent, as it can be easily verified. $\quad \square$

\noindent
Hence, the system of equations (\ref{vnl2}) is well posed, whatever the functions $\alpha_1$, $\dots$, $\alpha_\nu$ are; the already explicit version (\ref{constrexpl}) does not require any further condition.

\begin{rem}
In support of the form (\ref{ne}) of the equations of motion (in comparison with the more common equations in Lagrangian form, presented later on) is that they 
are promptly attainable, whenever the coefficients (\ref{gxi}) are simple: this is the case, for istance, of cartesian coordinates: if no holonomic constraint is present and $q_1, \dots, q_n$ are the cartesian coordinates ${\bf X}$ of the $N$ points, $n=3N$, formulae (\ref{coeffcdef}) are con\-si\-de\-ra\-bly simplified: actually, in (\ref{gxi}) the only non zero coefficients are $g_{i,i}$ containing 
(in triplets) the masses $M_1$, $\dots$, $M_N$ and
\begin{equation}
\label{coeffcart}
\begin{array}{l}
C_i^i=g_{i,i}+\sum\limits_{\nu=1}^k g_{m+\nu,m+\nu}
\left(\dfrac{\partial \alpha_\nu}{\partial {\dot q}_i}\right)^2, \quad 
C_i^r=C_r^i=\sum\limits_{\nu=1}^k g_{m+\nu,m+\nu}
\dfrac{\partial \alpha_\nu}{\partial {\dot q}_i}\dfrac{\partial \alpha_\nu}{\partial {\dot q}_r},
\quad \textrm{for}\;\;i\not=r \\
D_i^{r,s}=0, \quad E_i^r =  \sum\limits_\nu^k g_{m+\nu,m+\nu}\dfrac{\partial \alpha_\nu}{\partial q_r}\dfrac{\partial \alpha_\mu}{\partial {\dot q}_i}, \quad
G_i= \sum\limits_{\nu=1}^k g_{m+\nu, m+\nu}  \sum\limits_{\mu=1}^k\alpha_\mu 
\dfrac{\partial \alpha_\nu}{\partial q_{m+\mu}}.
\end{array}
\end{equation}
\end{rem}

\begin{exe} (nonholonomic pendulum, \cite{benpend}).
Two points $P_1$ and $P_2$ of mass $M_1$ and $M_2$ respectively are constrained on a $(x,y)$--plane, is such a way as to verify the following nonholonomic condition:
$$
{\dot y}_2(x_1{\dot x}_1+{\dot y}_1{\dot y}_1)-{\dot y}_1(x_2{\dot x}_2+y_2{\dot y}_2)=0
$$
where $(x_i, y_i)$ are the cartesian coordinates of $P_i$, $i=1,2$, 
Graphically, this means that the lines orthogonal to the velocities ${\dot P}_1$, ${\dot P}_2$ intersect in a point of the $y$--axis. We set $(q_1, q_2, q_3, q_4)=(x_1, y_1, y_2, x_2)$ and we put the kinematic constraint in the form (\ref{constrexpl}) by writing
$$
{\dot q}_4 = \dfrac{{\dot q}_3}{q_4} \left(q_2-q_3+\dfrac{q_1{\dot q}_1}{{\dot q}_2}\right) = \alpha_1 (q_1, q_2, q_3, q_4, {\dot q}_1, {\dot q}_2, {\dot q}_3).
$$
Setting 
$\Phi = \frac{q_1 {\dot q}_3}{q_4 {\dot q}_2}$, 
$\Psi =\frac{1}{q_4} \left(q_2-q_3+\frac{q_1{\dot q}_1}{{\dot q}_2}\right)$, 
coefficients (\ref{coeffcart}) are ($i$ denotes the row, $r$ the column)

$$
C_{i=1,2,3}^{r=1,2,3}=
\left(
\begin{array}{ccc}
M_1 + M_2 \Phi^2 & 
-M_2\frac{{\dot q}_1}{{\dot q}_2}
\Phi^2 & 
M_2 \Phi \Psi \\
&M_1+M_2\left(\frac{{\dot q}_1}{{\dot q}_2}\right)^2\Phi^2 & 
-M_2 \frac{{\dot q}_1}{{\dot q}_2}
\Phi \Psi
\\
& & M_2+M_2\Psi^2
\end{array}
\right)
$$

$$
E_{i=1,2,3}^{r=1,2,3}=m_2 \Phi 
\left(
\begin{array}{ccc}
\frac{{\dot q}_1}{q_1}\Phi & 
\frac{{\dot q}_3}{q_4} & 
-\frac{{\dot q}_3}{q_4}\\
-\frac{{\dot q}_1}{q_1}\frac{{\dot q}_1}{{\dot q}_2}\Phi  & -\frac{{\dot q}_1}{q_1}\Phi &
\frac{{\dot q}_1}{q_1}\Phi
 \\
\frac{{\dot q}_1}{q_1} \Psi &\frac{{\dot q}_2}{q_1} \Psi & -\frac{{\dot q}_2}{q_1} \Psi 
\end{array}
\right)
$$
$$
G=-M_2\frac{{\dot q}_3^2}{q_4} \Psi^2\left(\Phi, \frac{{\dot q}_1}{{\dot q}_2}\Phi, \Psi\right).
$$
After specifying the dynamics (spontaneous motion, gravitational field, ...) in a way that ${\cal F}^{(q_i)}$ can be written, the equations of motion (\ref{vnl2}) are readily achieved.
\end{exe}

\noindent
A different and more ordinary way to present the equations of motion for constrained system refers to the kinetic energy

\begin{equation}
\label{encin}
T(q_1,\dots, q_n, {\dot q}_1,\dots, {\dot q}_n,t)=
\frac{1}{2}{\mathcal Q}\cdot {\dot {\bf X}}, \qquad {\mathcal Q}={\dot {\bf X}}^{(M)}
\end{equation}
and to its restriction due to (\ref{constrexpl})
\begin{equation}
\label{trid}
T^*(q_1,\dots, q_n, {\dot q}_1, \dots, {\dot q}_m, t)=T(q_1, \dots, q_n, {\dot q}_1, \dots, {\dot q}_m, \alpha_1(\cdot), \dots, \alpha_k(\cdot), t)
\end{equation}
where each $\alpha_j(\cdot)$, $j=1,\dots, k$ depends on $q_1,\dots, q_n, {\dot q}_1, \dots, {\dot q}_m$ and $t$.
In that case, the equations of motion assume the form 

\begin{equation}
\label{vnl}
	\dfrac{d}{dt}\dfrac{\partial T^*}{\partial {\dot q}_i}-\dfrac{\partial T^*}{\partial q_i}
	-\sum\limits_{\nu=1}^k\dfrac{\partial T^*}{\partial q_{m+\nu}}\dfrac{\partial \alpha_\nu}{\partial {\dot q_i}}
	-\sum\limits_{\nu=1}^k  B_{i}^\nu \dfrac{\partial T}{\partial {\dot q}_{m+\nu}}
	={\cal F}^{(q_i)}+\sum\limits_{\nu=1}^k \dfrac{\partial \alpha_\nu}{\partial {\dot q}_i} {\cal F}^{(q_{m+\nu})}, 
	\qquad i=1,\dots, m
	\end{equation}
joined with (\ref{constrexpl}), where in $\dfrac{\partial T}{\partial {\dot q}_{m+\nu}}$ the variables ${\dot q}_{m+1}$, $\dots$, ${\dot q}_n$ are expressed in terms of $(q_1,\dots, q_n, {\dot q}_1, \dots, {\dot q}_m,t)$ by using (\ref{constrexpl}) and the coefficients $B_{i}^{\nu}(q_1,\dots, q_n, {\dot q}_1, \dots, {\dot q}_m,t)$ are defined by
	\begin{equation}
	\label{b}
	B_{i}^{\nu}
	=\sum\limits_{r=1}^m
	\left( \dfrac{\partial^2 \alpha_\nu}{\partial {\dot q}_i \partial q_r}{\dot q}_r +
	\dfrac{\partial^2 \alpha_\nu}{\partial {\dot q}_i \partial {\dot q}_r}{\q2dot^{..}}_r\right)
	-\dfrac{\partial \alpha_\nu}{\partial q_i}
	+\sum\limits_{\mu=1}^k\left( \dfrac{\partial^2 \alpha_\nu}{\partial {\dot q}_i \partial q_{m+\mu}}\alpha_\mu-
	\dfrac{\partial \alpha_\mu}{\partial {\dot q}_i}\dfrac{\partial \alpha_\nu}{\partial q_{m+\mu}}\right)+
	\dfrac{\partial^2 \alpha_\nu}{\partial {\dot q}_i \partial t}.
	\end{equation}
In order to achieve (\ref{vnl}), the key point is the well known relation 
${\dot {\mathcal Q}}\cdot \dfrac{\partial {\bf X}}{\partial q_i}=\dfrac{d}{dt}\dfrac{\partial T}{\partial {\dot q}_i}-\dfrac{\partial T}{\partial q_i}$. We remark that an alternative expression for 
$B_i^\nu$ is
\begin{equation}
\label{b2}B_i^\nu= \dfrac{d}{dt}\left( \dfrac{\partial \alpha_\nu}{\partial {\dot q}_i}\right)- 
\dfrac{\partial \alpha_\nu}{\partial q_i}-\sum\limits_{\mu=1}^k 
\dfrac{\partial \alpha_\mu}{\partial {\dot q}_i}
\dfrac{\partial \alpha_\nu}{\partial q_{m+\mu}}.
\end{equation}

\subsection{Some remarks}

\noindent
Equations (\ref{ne}) or (\ref{vnl}) hold for a system settled by $n$ parameters $q_1$, $\dots$, $q_n$ and undergoing the $m<n$ nonlinear nonholonomic constraints (\ref{constrexpl}); it is worth it to dwell upon some 
significant points and remarks.

\begin{description}

\item[$\bullet$] Equations (\ref{vnl2}) trace the Gibbs--Appell equations, since  
the left side of (\ref{vnl2}) corresponds to the calculus $\dfrac{\partial S}{\partial {\q2dot^{..}}_i}$, where $S=\frac{1}{2} {\dot {\cal Q}}\cdot \bfx2dot^{..}$ is the acceleration energy (Gibbs--Appell function); actually $\dfrac{\partial S}{\partial {\q2dot^{..}}_i}=
{\dot {\cal Q}}\cdot \dfrac{\partial \bfx2dot^{..}}{\partial {\q2dot^{..}}_i}={\dot {\cal Q}}\cdot
(\frac{\partial {\bf X}}{\partial q_i}+
\sum\limits_{j=1}^{k}\frac{\partial \alpha_j}{\partial {\dot q}_i}\frac{\partial {\bf X}}{\partial q_{m+j}})$, which sends back to (\ref{hatx}). 

\item[$\bullet$] On the other hand, equations (\ref{vnl2}) extend to the nonlinear case the Voronec equations (appeared in \cite{voronec}  and discussed  in the important monography \cite{neimark}) for the linear nonholonomic constraints (\ref{constrexpllin}); in the latter case the terms 
$\dfrac{\partial \alpha_\nu}{\partial {\dot q}_i}$ in (\ref{vnl2}) are simply $\alpha_{\nu,i}$ and
the equations of motion (\ref{vnl}) are, for each $i=1,\dots, m$, 
\begin{equation}
\label{voronec}
\dfrac{d}{dt}\dfrac{\partial T^*}{\partial {\dot q}_i}-\dfrac{\partial T^*}{\partial q_i}
-\sum\limits_{\nu=1}^k\alpha_{\nu,i}\dfrac{\partial T^*}{\partial q_{m+\nu}}
-\sum\limits_{\nu=1}^k\sum\limits_{r=1}^m 
\beta_{ir}^\nu{\dot q}_r\dfrac{\partial T}{\partial {\dot q}_{m+\nu}}=
{\mathcal F}^{(q_i)}+\sum\limits_{\nu=1}^{k} \alpha_{\nu,i} {\cal F}^{(q_{m+\nu})}
\end{equation}
where the coefficients (\ref{b}) reduce to 
\begin{equation}
\label{blin}
\begin{array}{l}
B_i^\nu=\sum\limits_{r=1}^m 
\beta_{ir}^\nu(q_1,\dots, q_n, t){\dot q}_r+\dfrac{\partial \alpha_{\nu,i}}{\partial t}, \\
\beta_{ir}^\nu=
\dfrac{\partial \alpha_{\nu, i}}{\partial q_r}-
\dfrac{\partial \alpha_{\nu,r}}{\partial q_i}+
\sum\limits_{\mu=1}^k\left(
\dfrac{\partial \alpha_{\nu, i}}{\partial q_{m+\mu}}\alpha_{\mu,r}-
\dfrac{\partial \alpha_{\nu, r}}{\partial q_{m+\mu}}\alpha_{\mu,i}
\right).
\end{array}
\end{equation}
Although the explicit dependence of $\alpha_{\nu,i}$ on $t$ is absent in \cite{neimark}, the widening to the rheonomic case is trivial.

\item[$\bullet$] 
Equations (\ref{vnl}) correspond to the ones derived in \cite{zek2}, as the most general form of equations of motion in Poincar\'e--Chetaev variables extended to nonlinear nonholonomic systems; the Voronec's 
equations (\ref{vnl}) are the same as the Voronec's equations pointed out in \cite{zek2} as the special case of Poincar\'e's kinematic parameters chosen as the real generalized velocities. 
Also the geometric approach for nonholonomic machanical systems (Lagrangian systems on fibered manifolds) performed in \cite{swa} leads to the same equations of motion as (\ref{vnl}). 

\item[$\bullet$] Concerning the dependence of (\ref{constrexpl}) on the variables, a special case is 
\begin{equation}
\begin{array}{ll}
\label{chap}
\alpha_\nu=\alpha_\nu(q_1, \dots, q_m, {\dot q}_1, \dots, {\dot q}_m,t)& \textrm{for each}\; \nu=1,\dots, k,\\
\\\
T=T(q_1, \dots, q_m, {\dot q}_1, \dots, {\dot q}_n, t), & \mathcal{F}_i=\mathcal{F}_i(q_1, \dots, q_m, {\dot q}_1, \dots, {\dot q}_n, t)
\end{array}
\end{equation}
that is the coordinates $q_{m+1}, \dots, q_n$ corresponding to the dependent velocities do not occur; 
we may refer to these systems as nonlinear $\check{\rm C}$aplygin systems. In this case, system (\ref{vnl}) reduces to (see also (\ref{b2}))

\begin{equation}
\label{chapeq}
\begin{array}{l}
\dfrac{d}{dt}\dfrac{\partial T^*}{\partial {\dot q}_i}-\dfrac{\partial T^*}{\partial q_i}
-\underbrace{\sum\limits_{\nu=1}^k  
\left( \sum\limits_{r=1}^m
\left( \dfrac{\partial^2 \alpha_\nu}{\partial {\dot q}_i \partial q_r}{\dot q}_r +
\dfrac{\partial^2 \alpha_\nu}{\partial {\dot q}_i \partial {\dot q}_r}{\q2dot^{..}}_r\right)
-\dfrac{\partial \alpha_\nu}{\partial q_i}
+\dfrac{\partial^2 \alpha_\nu}{\partial {\dot q}_i \partial t}\right)}_{\frac{d}{dt}
\frac{\partial \alpha_\nu}{\partial {\dot q}_i}-\alpha_\nu}
 \dfrac{\partial T}{\partial {\dot q}_{m+\nu}}\\
={\cal F}^{(q_i)}+\sum\limits_{\nu=1}^k \dfrac{\partial \alpha_\nu}{\partial {\dot q}_i} {\cal F}^{(q_{m+\nu})}.
\end{array}
\end{equation}
The remarkable advantage of the differential
system (\ref{chapeq}) is that it contains only the unknown functions $q_1$, $\dots$, $q_m$ and it is disentangled from the constraints equations (\ref{constrexpllin}).

\noindent
Within assumption (\ref{chap}), 
the linear stationary case 
\begin{equation}
\begin{array}{ll}
\label{chaplin}
\alpha_{\nu,j}=\alpha_{\nu,j} (q_1, \dots, q_m)& \textrm{for each}\; \nu=1,\dots, k\;\textrm{and}\;j=1, \dots,m\\
\\\
T=T(q_1, \dots, q_m, {\dot q}_1, \dots, {\dot q}_n), & \mathcal{F}_i=\mathcal{F}_i(q_1, \dots, q_m, {\dot q}_1, \dots, {\dot q}_n)
\end{array}
\end{equation}
leads to the $\check{\rm C}$aplygin's equations (see (\cite{neimark}))
\begin{equation}
\label{capligineq}
\dfrac{d}{dt}\dfrac{\partial T^*}{\partial {\dot q}_i}-\dfrac{\partial T^*}{\partial q_i}
-\sum\limits_{\nu=1}^k\sum\limits_{r=1}^m 
\left(\dfrac{\partial \alpha_{\nu, i}}{\partial q_r}-
\dfrac{\partial \alpha_{\nu,r}}{\partial q_i}\right)
{\dot q}_r\dfrac{\partial T}{\partial {\dot q}_{m+\nu}}=
{\cal F}^{(q_i)}+\sum\limits_{j=1}^{k} \alpha_{j,i} {\cal F}^{(q_{m+j})}
\quad i=1,\dots, m.
\end{equation}

\item[$\bullet$]
Whenever the kinetic energy (\ref{encin}) is, consistently with the usual mechanical systems,
$$
T=\frac{1}{2}
\sum\limits_{i,j=1}^n g_{i,j}{\dot q}_i {\dot q}_j +\sum\limits_{i=1}^n b_i {\dot q}_i+c, 
\quad
b_i(q_1,\dots, q_n,t)=\dfrac{\partial {\bf X}^{(\mathsf{M})}}{\partial q_i}\cdot \dfrac{\partial {\bf X}}{\partial t}, \;
c(q_1,\dots, q_n,t)=\frac{1}{2}
\dfrac{\partial {\bf X}^{(\mathsf{M})}}{\partial t}\cdot 
\dfrac{\partial {\bf X}}{\partial t}
$$
($g_{i,j}$ are defined in (\ref{gxi})) so that (\ref{trid}) writes
\begin{equation}
\label{tmobile}
T^*=
\frac{1}{2}\left( \sum\limits_{r,s=1}^m g_{r,s}{\dot q}_r {\dot q}_s+
\sum\limits_{\nu, \mu=1}^k g_{m+\nu,m+\mu} \alpha_\nu \alpha_\mu \right)+
\sum\limits_{r=1}^m \sum\limits_{\nu=1}^k g_{r,m+\nu}{\dot q}_r \alpha_\nu+
\sum\limits_{r=1}^m b_r {\dot q}_r +\sum\limits_{\nu=1}^k b_{m+\nu}\alpha_\nu+c
\end{equation}
then rearranging the terms in (\ref{vnl}) one can easily check (see \cite{tal})) that all the terms of $-\sum\limits_{\nu=1}^k \dfrac{\partial T}{\partial {\dot q}_{m+\nu}} B_{i}^\nu$, for each $i=1,\dots, m$, cancel out with part of the addends of $\dfrac{d}{dt}\left(\dfrac{\partial T^*}{\partial {\dot q}_i}\right)$,  
	of $-\dfrac{\partial T^*}{\partial q_i}$ and of
	$-\sum\limits_{\nu=1}^k\dfrac{\partial T^*}{\partial q_{m+\nu}}\dfrac{\partial \alpha_\nu}{\partial {\dot q_i}}$. 
The remaining terms of (\ref{vnl}) coincide precisely with the terms of (\ref{vnl2}).
\end{description}

\noindent
In order to clarify the latter statement, we illustrate the following

\begin{exe}
On a plane $(x,y)$ we consider two points $P_1\equiv(x_1, y_1)$ and $P_2\equiv (x_2, y_2)$ whose velocities are perpendicular: 
\begin{equation}
\label{velperp}
{\dot x}_1 {\dot x}_2+{\dot y}_1{\dot y}_2=0.
\end{equation}
We set ${\bf X}(q_1, q_2, q_3, q_4)$, $q_1=x_1$, $q_2=y_2$, $q_3=x_2$, $q_4=y_2$, so that (\ref{constrexpl}) is 
$$
{\dot q}_4=\alpha_1({\dot q}_1, {\dot q}_2, {\dot q}_3)=-\frac{{\dot q}_1 {\dot q}_3}{{\dot q}_2}.
$$ 
The system runs into the case (\ref{coeffcart}) with $g_{1,1}=g_{2,2}=M_1$, $g_{3,3}=g_{4,4}=M_2$ 
and the coeffcients are 
$$
\begin{array}{lll}
C_1^1= M_1 + M_2 ({\dot q}_3/{\dot q}_2)^2 & 
C_2^2= M_1 + M_2 ({\dot q}_1{\dot q}_3/{\dot q}_2^2)^2 &
C_3^3= M_1 +  M_2 ({\dot q}_1/{\dot q}_2)^2\\
C_1^2=C_2^1= - M_2 {\dot q}_1 {\dot q}_3^2/{\dot q}_2^3&
C_1^3=C_3^1= M_2 {\dot q}_1 {\dot q}_3^2/{\dot q}_2^2&
C_2^3=C_3^2 = -M_2{\dot q}_1^2{\dot q}_3/{\dot q}_2^3  \\
D_i^{r,s}=E_i^r=G_i=0 & & 
\end{array}
$$
so that (\ref{vnl2}) are (let us leave the forces terms unspecified)
\begin{equation}
\label{vnl2velperp}
\begin{array}{l}
(M_1 + M_2 ({\dot q}_3/{\dot q}_2)^2) {\q2dot^{..}}_1
-  M_2({\dot q}_1 {\dot q}_3/{\dot q}_2^3)({\dot q}_3{\q2dot^{..}}_2- {\dot q}_2{\q2dot^{..}}_3)
={\cal F}^{(q_1)}-{\cal F}^{(q_4)} {\dot q}_3/ {\dot q}_2\\
\\
(M_1 +  M_2 ({\dot q}_1{\dot q}_3/{\dot q}_2^2)^2){\q2dot^{..}}_2-
M_2 ({\dot q}_1 {\dot q}_3/{\dot q}_2^3) ({\dot q}_3  {\q2dot^{..}}_1-
{\dot q}_1  {\q2dot^{..}}_3)
={\cal F}^{(q_1)}+{\cal F}^{(q_4)}{\dot q}_1{\dot q}_3/ {\dot q}_2^2 \\
\\
(M_1 +  M_2 ({\dot q}_1/{\dot q}_2)^2) {\q2dot^{..}}_3+
M_2({\dot q}_1 {\dot q}_3/{\dot q}_2^3)({\dot q}_2{\q2dot^{..}}_1-{\dot q}_1{\q2dot^{..}}_2)
={\cal F}^{(q_1)}-{\cal F}^{(q_4)} {\dot q}_1/ {\dot q}_2
\end{array}
\end{equation}
On the other hand, by writing (\ref{trid}) as $T^*=\frac{1}{2}M_1({\dot q}_1^2
+{\dot q}_2^2)+\frac{1}{2}M_2{\dot q}_3^2(1+({\dot q}_1/{\dot q}_2)^2)$ and calculating (\ref{b}) as
$$
B_1^1=({\dot q}_3{\q2dot^{..}}_2-{\dot q}_2{\q2dot^{..}}_3)/{\dot q}_2^2, \quad 
B_2^1=({\dot q}_1{\q2dot^{..}}_3+{\dot q}_3{\q2dot^{..}}_1)/{\dot q}_2^2
-({\dot q}_1 {\dot q}_3/{\dot q}_2^3){\q2dot^{..}}_2, \quad
B_3^1=(-{\dot q}_2{\q2dot^{..}}_1+{\dot q}_1{\q2dot^{..}}_2)/{\dot q}_2^2
$$
the calculation of (\ref{vnl}) leads to
\begin{equation}
\label{vnlvelperp}
\left\{
\begin{array}{l}
\dfrac{d}{dt}[{\dot q}_1(M_1+M_2 ({\dot q}_3/{\dot q}_2)^2)]
+M_2({\dot q}_1{\dot q}_3/{\dot q}_2^3)({\dot q}_3{\q2dot^{..}}_2- {\dot q}_2{\q2dot^{..}}_3)=
{\cal F}^{(q_1)}- {\cal F}^{(q_4)} {\dot q}_3/ {\dot q}_2\\
\\
\dfrac{d}{dt}[{\dot q}_2(M_1-M_2 ({\dot q}_1{\dot q}_3/{\dot q}_2^2)^2)]
+M_2({\dot q}_1{\dot q}_3/{\dot q}_2^3)({\dot q}_3{\q2dot^{..}}_1+{\dot q}_1{\q2dot^{..}}_3
-({\dot q}_1{\dot q}_3/{\dot q}_2){\q2dot^{..}}_2)
={\cal F}^{(q_2)}+{\cal F}^{(q_4)}{\dot q}_1{\dot q}_3/ {\dot q}_2^2\\
\\
\dfrac{d}{dt}[M_2{\dot q}_3(1+({\dot q}_1/{\dot q}_2)^2)]
-M_2({\dot q}_1{\dot q}_3/{\dot q}_2^3)({\dot q}_2{\q2dot^{..}}_1- {\dot q}_1{\q2dot^{..}}_2)
={\cal F}^{(q_3)}-{\cal F}^{(q_4)} {\dot q}_1/ {\dot q}_2.
\end{array}
\right.
\end{equation}
According to the assertion claimed after (\ref{tmobile}), all the terms 
$-\dfrac{\partial T}{\partial {\dot q}_4} B_{i}^1$, $i=1,2,3$, which are the terms beginning with $+m_2$ in each of (\ref{vnlvelperp}), cancel out with opposite terms of the explicit calculation of $d/dt$. 
The rest of the terms are those appearing in (\ref{vnl2velperp}).

\end{exe}

\begin{exe} 
A different example, discussed in \cite{virga}, draws attention to the fact that equations (\ref{vnl}) are still valid even though the function $T^*$ defined in (\ref{tmobile}) degenerates w.~r.~t.~the restricted variables ${\dot q}_1$, $\dots$, ${\dot q}_m$: let us consider one point of mass $M$ of cartesian coordinates $(x,y,z)$ and whose velocity is constant in module: ${\dot x}^2+{\dot y}^2+{\dot z}^2=C^2\not = 0$. In terms of $(q_1, q_2, q_3)=(x,y,z)$ we write (\ref{constrexpl})
as ${\dot q}_3=\pm \sqrt{C^2-{\dot q}_1^2 - {\dot q}_2^2}=\alpha_1({\dot q}_1, {\dot q}_2)$, where the sign derives from the initial conditions (this example presents $m=2$, $k=1$). 
Concerning (\ref{vnl2}), we use (\ref{coeffcart}) with $g_{1,1}=g_{2,2}=g_{3,3}=m$, $g_{i,j}=0$ for $i\not =j$, so that
$$
C_i^i= M\left(
1+\frac{M{\dot q}_i^2}{C^2-{\dot q}_1^2-{\dot q}_2^2}\right), \quad i=1,2, \qquad  
C_2^1=C_1^2=
\frac{M {\dot q}_1 {\dot q}_2}{C^2- {\dot q}_1^2 - {\dot q}_2^2}
$$
so that equations (\ref{vnl}) are
\begin{equation}
\label{velcost}
\left\{
\begin{array}{l}
M \frac{C^2-{\dot q}_2^2}{C^2-{\dot q}_1^2 -{\dot q_2^2}}{\q2dot^{..}}_1
+ M \frac{{\dot q}_1 {\dot q}_2}{C^2-{\dot q}_1^2 -{\dot q_2^2}}{\q2dot^{..}}_2
= {\cal F}^{(q_1)}\mp {\cal F}^{(q_3)} 
\frac{{\dot q}_1}{\sqrt{C^2-{\dot q}_1^2 - {\dot q}_2^2}}, \\ 
M \frac{{\dot q}_1 {\dot q}_2}{C^2-{\dot q}_1^2 -{\dot q}_2^2}{\q2dot^{..}}_1+
M \frac{C^2-{\dot q}_1^2}{C^2-{\dot q}_1^2 -{\dot q}_2^2}{\q2dot^{..}}_2
= {\cal F}^{(q_1)}\mp {\cal F}^{(q_3)} 
\frac{{\dot q}_2}{\sqrt{C^2-{\dot q}_1^2 - {\dot q}_2^2}},	
\end{array}
\right.
\end{equation}
On the other hand, (\ref{tmobile}) is $T^*=\frac{1}{2}mC^2$, hence the only contributions in the left hand side of (\ref{vnl}) are 
$$
-B^1_i\dfrac{\partial T}{\partial {\dot q}_3}=
-m\left(\dfrac{\partial^2 \alpha_1}{\partial {\dot q}_i \partial {\dot q}_1}{\q2dot^{..}}_1 
+\dfrac{\partial^2 \alpha_1}{\partial {\dot q}_i \partial {\dot q}_2}{\q2dot^{..}}_2\right){\dot q}_3({\dot q}_1, {\dot q}_2), \quad i=1, 2
$$
Calculating the second derivatives, one easily finds exactly equations (\ref{velcost}).

\end{exe}

\section{Energy balance}

\noindent
Let us assume that the active forces depend only on ${\bf X}$ and $t$ and come from a potential ${\cal U}$:
\begin{equation}
\label{gradu}
{\mathcal F}=\nabla_{\bf X}{\cal U}({\bf X}, t)
\end{equation} 
so that the restriction $U(q_1, \dots, q_n,t)={\cal U}({\bf X}(q_1, \dots, q_n,t),t)$ to the configuration manifold provides the generalized forces as follows:
\begin{equation}
\label{potu}
{\mathcal F}^{(q_i)}=\dfrac{\partial U}{\partial q_i}, \quad i=1, \dots, n\quad 
\end{equation}

\noindent
It is known that in case of holonomic systems ${\bf X}(q_1, \dots, q_\ell,t)$ (that is removing (\ref{constr})) the equations of motion $\dfrac{d}{dt}\dfrac{\partial {\cal L}}{\partial {\dot q}_i}-\dfrac{\partial {\cal L}}{\partial q_i}=0$, $i=1, \dots, n$, where ${\cal L}=T+U$ (see (\ref{encin})) entail the energy balance $\dfrac{d}{dt}\left( \sum\limits_{i=1}^n {\dot q}_i \dfrac{\partial {\cal L}}{\partial {\dot q}_i}-{\cal L}\right)=-\dfrac{\partial {\cal L}}{\partial t}$, which supplies the conservation of the quantity in brackets, whenever ${\cal L}$ does not depend on $t$ explicitly.

\noindent
Now, if the constraints (\ref{constr}) are present, recalling $T^*$ defined in (\ref{trid}) we set
\begin{equation}
\label{lstar}
{\cal L}^*(q_1, \dots, q_n, {\dot q}_1,\dots, {\dot q}_m, t)=T^*(q_1, \dots, q_n, {\dot q}_1,\dots, {\dot q}_m, t)+U(q_1, \dots, q_n,t)
\end{equation}
as the Lagrangian in terms of the independent velocities. The following Proposition generalizes the just mentioned balance of holonomic systems.

\begin{prop} 
The equations of motion (\ref{vnl}) entail
\begin{equation}
\label{bilen}
\dfrac{d}{dt}\left( \sum\limits_{i=1}^m {\dot q}_i \dfrac{\partial {\cal L}^*}{\partial {\dot q}_i}-{\cal L}^*\right)-
	\sum\limits_{\nu=1}^k ({\overline \alpha}_\nu-\alpha_\nu)\dfrac{\partial {\cal L}^*}{\partial q_{m+\nu}}-
	\sum\limits_{\nu=1}^k {\overline B}_\nu \dfrac{\partial T}{\partial {\dot q}_{m+\nu}}=
	-\dfrac{\partial {\cal L}^*}{\partial t}
\end{equation}
where
\begin{eqnarray}
\label{baralpha}
{\overline \alpha}_\nu (q_1, \dots, q_n, {\dot q}_1, \dots, {\dot q}_m,t)&=&
\sum\limits_{i=1}^m 
{\dot q}_i\dfrac{\partial \alpha_\nu}{\partial {\dot q}_i}, \\
{\overline B}_\nu (q_1, \dots, q_n, {\dot q}_1, \dots, {\dot q}_m,t)
\label{barbnu}
&=&\sum\limits_{r=1}^m {\dot q}_r \left( \dfrac{\partial {\overline \alpha}_\nu}{\partial q_r}-
\dfrac{\partial \alpha_\nu}{\partial q_r}\right)
+\sum\limits_{\mu=1}^k \left( \alpha_\mu \dfrac{\partial {\overline \alpha}_\nu}{\partial q_{m+\mu}}-
{\overline \alpha}_\mu \dfrac{\partial \alpha_\nu}{\partial q_{m+\mu}}\right)\\
\nonumber
&+&\sum\limits_{r=1}^m {\q2dot^{..}}_r \left( 
\dfrac{\partial {\overline \alpha}_\nu }{\partial {\dot q}_r} -
\dfrac{\partial \alpha_\nu }{\partial {\dot q}_r}\right)
+\dfrac{\partial {\overline \alpha}_\nu}{\partial t}\\
\nonumber
&=&\dfrac{d}{dt}({\overline \alpha}_\nu-\alpha_\nu)-\sum\limits_{\mu=1}^k
\dfrac{\partial \alpha_\nu}{\partial q_{m+\mu}}
({\overline \alpha}_\mu -\alpha_\mu)+\dfrac{\partial \alpha_\nu}{\partial t}.
\end{eqnarray}
\end{prop}

\noindent
{\bf Proof}. Basing on the formula 
\begin{equation}
\label{f}
\dfrac{d}{dt}F(q_1, \dots, q_n, {\dot q}_1, \dots, {\dot q}_m,t)=
\sum\limits_{i=1}^m \dfrac{\partial F}{\partial q_i}{\dot q}_i+\sum\limits_{\nu=1}^k 
\dfrac{\partial F}{\partial q_{m+\nu}}\alpha_\nu+
\sum\limits_{i=1}^m \dfrac{\partial F}{\partial {\dot q}_i} {\q2dot^{..}}_i+
\dfrac{\partial F}{\partial t}
\end{equation}
implemented with $F=T^*$ 
one finds 
$$ 
\sum\limits_{i=1}^m
{\dot q}_i	
\left(
\dfrac{d}{dt}\dfrac{\partial T^*}{\partial {\dot q}_i}-\dfrac{\partial T^*}{\partial q_i}
-\sum\limits_{\nu=1}^k\dfrac{\partial T^*}{\partial q_{m+\nu}}\dfrac{\partial \alpha_\nu}{\partial {\dot q_i}}
\right)
=\dfrac{d}{dt}\left( \sum\limits_{i=1}^m {\dot q}_i \dfrac{\partial T^*}{\partial {\dot q}_i}-T^*\right)
+\dfrac{\partial T^*}{\partial t}
-\sum\limits_{\nu=1}^k \dfrac{\partial T^*}{\partial q_{m+\nu}}
\left(\sum\limits_{i=1}^m
{\dot q}_i\dfrac{\partial \alpha_\nu}{\partial {\dot q}_i} -\alpha_\nu\right)
$$
\begin{equation}
\label{ultimaa}
=\dfrac{d}{dt}\left( \sum\limits_{i=1}^m {\dot q}_i \dfrac{\partial T^*}{\partial {\dot q}_i}-T^*\right)
+\dfrac{\partial T^*}{\partial t}
-\sum\limits_{\nu=1}^k \dfrac{\partial T^*}{\partial q_{m+\nu}}
\left({\overline \alpha}_\nu -\alpha_\nu\right).
\end{equation}
Furthermore, recalling (\ref{b2}) it is
\begin{eqnarray}
\nonumber
\sum\limits_{i=1}^m B_i^\nu {\dot q}_i&=&
\dfrac{d}{dt}\left(\sum\limits_{i=1}^m {\dot q}_i \dfrac{\partial \alpha_\nu}{\partial {\dot q}_i}\right)- 
\sum\limits_{i=1}^m\left(
{\dot q}_i \dfrac{\partial \alpha_\nu}{\partial q_i}+
{\q2dot^{..}}_i \dfrac{\partial \alpha_\nu}{\partial {\dot q}_i}\right)
-\sum\limits_{\mu=1}^k 
\sum\limits_{i=1}^m
{\dot q}_i\dfrac{\partial \alpha_\mu}{\partial {\dot q}_i}
\dfrac{\partial \alpha_\nu}{\partial q_{m+\mu}} \\
&=&
\nonumber
\dfrac {d {\overline \alpha}_\nu}{dt}-
\sum\limits_{i=1}^m \left(
{\dot q}_i \dfrac{\partial \alpha_\nu}{\partial q_i}
+{\q2dot^{..}}_i \dfrac{\partial \alpha_\nu}{\partial {\dot q}_i}\right)
-\sum\limits_{\mu=1}^k 
{\overline \alpha}_\mu
\dfrac{\partial \alpha_\nu}{\partial q_{m+\mu}}
\end{eqnarray}
\begin{equation}
\label{ultimab}
=
\sum\limits_{r=1}^m {\dot q}_r \left( \dfrac{\partial {\overline \alpha}_\nu}{\partial q_r}
- \dfrac{\partial \alpha_\nu}{\partial q_r}\right)
+\sum\limits_{\mu=1}^k \left( \alpha_\mu 
\dfrac{\partial  {\overline \alpha}_\nu}{\partial q_{m+\mu}}
- {\overline \alpha}_\mu \dfrac{\partial \alpha_\nu}{\partial q_{m+\mu}}\right)
+\sum\limits_{r=1}^m {\q2dot^{..}}_r \left( 
\dfrac{\partial {\overline \alpha}_\nu }{\partial {\dot q}_r} -
\dfrac{\partial \alpha_\nu }{\partial {\dot q}_r}\right)
+\dfrac{\partial {\overline \alpha}_\nu }{\partial t}.
\end{equation}
Concerning the forces, under assumption (\ref{potu}) and having in mind $\dfrac{dU}{dt}=
\sum\limits_{i=1}^m \dfrac{\partial U}{\partial q_i}{\dot q}_i+\sum\limits_{\nu=1}^k \dfrac{\partial U}{\partial q_{m+\nu}}\alpha_\nu+\dfrac{\partial U}{\partial t}$, 
we can write
\begin{equation}
\label{balu}
\sum\limits_{i=1}^m {\dot q}_i	
\left(
{\cal F}^{(q_i)}+\sum\limits_{\nu=1}^k \dfrac{\partial \alpha_\nu}{\partial {\dot q}_i} {\cal F}^{(q_{m+\nu})}
\right)=\dfrac{dU}{dt}-\dfrac{\partial U}{\partial t}+\sum\limits_{\nu=1}^k \dfrac{\partial U}{\partial q_{m+\nu}}\left({\overline \alpha}_\nu -\alpha_\nu\right).
\end{equation}
By virtue of (\ref{ultimaa}), (\ref{ultimab}) and (\ref{balu}), multiplying each of (\ref{vnl}) by ${\dot q}_i$ and summing up with respect to $i$ one gets the statement (\ref{bilen}). The second equality for ${\overline B}_\nu$ in (\ref{barbnu}) is obtained by applying (\ref{f}) with $F={\overline \alpha}_\nu - \alpha_\nu$. $\quad \square$

\begin{cor}
For a system verifying assumption (\ref{chap}) (nonlinear $\check{\rm C}$aplygin's systems) and assumption (\ref{potu}) for $i=1, \dots, m$, equation (\ref{bilen}) takes the simpler form
	\begin{equation}
	\label{bilenchap}
	\dfrac{d}{dt}\left( \sum\limits_{i=1}^m {\dot q}_i \dfrac{\partial {\cal L}^*}{\partial {\dot q}_i}-{\cal L}^*\right)-
	\sum\limits_{\nu=1}^k \left( 
\dfrac{d}{dt}({\overline \alpha}_\nu-\alpha_\nu)+\dfrac{\partial \alpha_\nu}{\partial t}		
	\right) \dfrac{\partial T}{\partial {\dot q}_{m+\nu}}=
	-\dfrac{\partial {\cal L}^*}{\partial t}
	\end{equation}
		
\end{cor}
Indeed, the terms containing $\dfrac{\partial {\cal L}^*}{\partial q_{m+\nu}}$ in (\ref{baralpha}) cancel out for each $\nu=1,\dots, k$.

\noindent
Whenever $T^*$ is the function (\ref{tmobile}), the energy of the system is
\begin{eqnarray}
\label{enexpl}
\sum\limits_{i=1}^m {\dot q}_i \dfrac{\partial {\cal L}^*}{\partial {\dot q}_i}-{\cal L}^*&=&
\frac{1}{2}\sum\limits_{r,s}^m 
g_{r,s}{\dot q}_r {\dot q}_s
+\sum\limits_{\nu, \mu=1}^k 
g_{m+\nu, m+\mu} ({\overline \alpha}_\nu \alpha_\mu - \frac{1}{2}\alpha_\nu \alpha_\mu)\\
\nonumber
&+&\sum\limits_{r=1}^m \sum\limits_{\nu=1}^k g_{r,m+\nu} {\overline \alpha}_\nu {\dot q}_r 
+\sum\limits_{\nu=1}^k b_{m+\nu} ({\overline \alpha}_\nu- \alpha_\nu)-c-U.
\end{eqnarray}
We refer to (\ref{enexpl}) as the energy of the system. 

\begin{exe}
Let us exert (\ref{bilen}) for the system of Example 1 (nonholonomic pendulum): 
assuming that the forces give rise to the potential $U(q_1, q_2, q_3, q_4)$, the function (\ref{lstar})
is, recalling that the only one kinematic constraint is ${\dot q}_4 =  \dfrac{{\dot q}_3}{q_4} \left(q_2-q_3+\dfrac{q_1{\dot q}_1}{{\dot q}_2}\right) = \alpha_1 (q_1, q_2, q_3, q_4, {\dot q}_1, {\dot q}_2, {\dot q}_3)$,
$$
{\cal L}^*=\dfrac{1}{2}M_1 ({\dot q}_1^2+{\dot q}_2^2)+\dfrac{1}{2}M_2 {\dot q}_3^2\left(1+ \dfrac{1}{q_4^2} \left(q_2-q_3+\dfrac{q_1 {\dot q}_1}{{\dot q}_2}\right)^2\right)+U(q_1, q_2, q_3, q_4).
$$
Calculating (\ref{baralpha}) for $\alpha_1$ one finds
$$
{\overline {\bar \alpha}_1}=\cancel{\frac{q_1{\dot q}_3}{q_4{\dot q}_2}{\dot q}_1} + 
\cancel{\frac{{\dot q}_3}{q_4}\left(
-\dfrac{q_1 {\dot q}_1}{{\dot q}_2^2}\right){\dot q}_2} 
+ \dfrac{1}{q_4} \left(q_2-q_3+\dfrac{q_1{\dot q}_1}{{\dot q}_2}\right){\dot q}_3 = \alpha_1
$$
so that even ${\overline B}_1 =0$ (see (\ref{barbnu})). Therefore, (\ref{bilen}) provides the constant of motion (\ref{enexpl}) (in this example $g_{1,1}=g_{2,2}=M_1$, $g_{3,3}=g_{4,4}=M_2$, $g_{i,j}=0$ for $i\not =j$)
$$
I(q_1, q_2, q_3, q_4, {\dot q}_1, {\dot q}_2, {\dot q}_3)=
\dfrac{1}{2}M_1 ({\dot q}_1^2+{\dot q}_2^2)+\dfrac{1}{2}M_2 
+\dfrac{1}{2}M_2 {\dot q}_3^2\left(1+ \dfrac{1}{q_4^2} \left(q_2-q_3+\dfrac{q_1 {\dot q}_1}{{\dot q}_2}\right)^2\right)-U(q_1, q_2, q_3, q_4).
$$
\end{exe}

\begin{exe}
$N$ material points $(P_i, M_i)$, $i=1,\dots, N$ lie on the $(x,y)$--plane and are constrained to move with parallel velocities; denoting by $(x_i, y_i)$, $i=1,2,3$, the cartesian coordinates of the points, the nonholonomic restrictions are
$$
{\dot x}_1{\dot y}_2-{\dot x}_2{\dot y}_1=0, \quad {\dot x}_1{\dot y}_3-{\dot x}_3{\dot y}_1=0 \quad \dots \quad {\dot x}_1{\dot y}_N-{\dot x}_N{\dot y}_N=0 
$$  
which can be exhibited in the explicit form (\ref{constrexpl}) by setting, for $n=2N$, $(q_1, q_2, \dots,  q_n)=
	(\overbrace{x_1, x_2, \dots, x_N, y_1}^{q_1, \dots, q_m}, \underbrace{y_2, \dots, y_N}_{q_{m+1}, \dots, q_n})$, and by writing 
	$$
	{\dot q}_{m+1}=\frac{{\dot q}_2}{{\dot q}_1}{{\dot q}_m}, \quad \dots \;\;
	{\dot q}_n=\frac{{\dot q}_{m-1}}{{\dot q}_1}{{\dot q}_m}
	$$
	that is, for $\nu=1, \dots, k$, ${\dot q}_{m+\nu}=\alpha_\nu=\frac{{\dot q}_{\nu+1}}{{\dot q}_1}{{\dot q}_m}$; in this example $m=N+1$, $k=N-1$. The coefficients (\ref{baralpha}) and (\ref{barbnu}) are, for each $\nu=1, \dots, m$, 
	$$
	{\overline \alpha}_\nu=\cancel{-\dfrac{{\dot q}_{\nu+1}}{{\dot q}_1^2}{\dot q}_m{\dot q}_1}+\cancel{\dfrac{{\dot q}_{\nu+1}}{{\dot q}_2}}+\dfrac{{\dot q}_m}{{\dot q}_1}{\dot q_{\nu+1}}=\alpha_\nu, \qquad {\overline B}_\nu =0.
	$$
	Therefore
	${\cal L}^*=\dfrac{1}{2}\sum\limits_{i=1}^{m-1} M_i{\dot q}_i^2 (1+{\dot q}_m^2/{\dot q}_1^2)$
	gives by means of (\ref{bilen}) the constant of motion (\ref{enexpl}): having in mind 
	$g_{1,1}=g_{m,m}=M_1$, $g_{2,2}=g_{m+1, m+1}=M_2$, $\dots$, $g_{m-1,m-1}=g_{n,n}=M_N$, elsewhere null, 
	the first integral is $\dfrac{1}{2}\sum\limits_{i=1}^{m-1} M_i{\dot q}_i^2 (1+{\dot q}_m^2/{\dot q}_1^2)-U(q_1, \dots, q_n)$, where the last term takes into account the possibile interactions. 
\end{exe}

\begin{exe}
We go back to Example 3 (one point with constant norm of the velocity), in order to write (\ref{bilen}) for (\ref{velcost}). Assuming that the active forces are connected to the potential $U=U(q_1, q_2, q_3)$, we have ${\cal L}^*= \frac{1}{2}MC^2+U$
and (\ref{bilen}) writes
$$-\dfrac{dU}{dt}-({\overline \alpha}_1 - \alpha_1)\dfrac{\partial U}{\partial q_3}-M
\alpha_1{\overline B}_1=0
$$ 
with $\alpha_1= \pm \sqrt{C^2-{\dot q}_1^2-{\dot q}_2^2}$,
${\overline \alpha}_1= \mp\frac{{\dot q}_1^2+{\dot q}_2^2}{\sqrt{C^2-{\dot q}_1^2-{\dot q}_2^2}}$,  
${\overline B}_1=\mp \frac{C^2}{(C^2-{\dot q}_1^2 -{\dot q}_2^2)^{3/2}}
({\dot q}_1{\q2dot^{..}}_1+ {\dot q}_2{\q2dot^{..}}_2)$. Hence we get
$$
-\dfrac{\partial U}{\partial q_1}{\dot q}_1 - \dfrac{\partial U}{\partial q_2}{\dot q}_2 
\pm \frac{{\dot q}_1^2+{\dot q}_2^2}{\sqrt{C^2-{\dot q}_1^2-{\dot q}_2^2}}\dfrac{\partial U}{\partial q_3}
=\dfrac{M}{2} \dfrac{C^2}{C^2-{\dot q}_1^2-{\dot q}_2^2}  \dfrac{d}{dt}({\dot q}_1^2+{\dot q}_2^2).
$$
\end{exe}

\begin{exe}
We can modify the previous example by requiring two material points to have the same norm of the velocity: 
${\dot P}_1^2= {\dot P}_2^2$. Calling $(x_i, y_i, z_i)$ the cartesian coordinate of $P_i$, $i=1,2$, and setting $(q_1, q_2, q_3, q_4, q_5, q_6)=(x_1, y_1, z_1, x_2, y_2, z_2)$, the form (\ref{constrexpl})
of the nonholonomic restriction is 
$$
{\dot q}_6 =\pm \sqrt{{\dot q}_1^2+ {\dot q}_2^2 +{\dot q}_3^2 - {\dot q}_4^2-{\dot q}_5^2}=
\alpha_1({\dot q}_1, {\dot q}_2, {\dot q}_3,{\dot q}_4, {\dot q}_5).
$$ 
Even in this case ${\overline \alpha}_1=\alpha_1$, ${\overline B}_1=0$ and (\ref{bilen})  entails the conservation of the quantity (\ref{enexpl}) $\frac{1}{2}(M_1+M_2)({\dot q}_1^2+{\dot q}_2^2 +{\dot q}_3^2) - U(q_1, q_2, q_3, q_4, q_5, q_6)$, whenever the forces allow the access to the potential $U$.  
	
\end{exe}

\section{Special classes of nonholonomic constraints}

\noindent
As it emerges from the Examples, the energy balance (\ref{bilen}) deserves a distinctive treatment whenever the constraint functions (\ref{constrexpl}) take a specific form. In particular, the circumstance ${\overline \alpha}_\nu=\alpha_\nu$ play the crucial role for the conservation of the energy of the system. Let us start from the following

\begin{lem}
	For a fixed $\nu$ from $1$ up to $k$, ${\overline \alpha}_\nu=\alpha_\nu$ if and only if $\alpha_\nu$ is a homogeneous function of degree $1$ w.~r.~t.~${\dot q}_1$, $\dots$, ${\dot q}_m$.
\end{lem}

\noindent
{\bf Proof}. We simply turn to the Euler's theorem: 
${\overline \alpha}_\nu=\sum\limits_{r=1}^m {{\dot q}_r}\dfrac{\partial \alpha_\nu}{\partial {\dot q}_r}=\alpha_\nu$ if and only if $\alpha_\nu$ is a homogeneous function of degree $1$ with respect to the variables ${\dot q}_1$, $\dots$, ${\dot q}_m$. $\quad\square$

\noindent
Let us now assume that each of the nonholonomic constraints verifies
\begin{equation}
\label{baralphaalpha}
{\overline \alpha}_\nu =\alpha_\nu \quad \textrm{for any}\;\;\nu=1, \dots, k.
\end{equation}
Then, the energy balance simplifies according to the following statement.

\begin{prop}
If (\ref{baralphaalpha}): holds, the energy balance (\ref{bilen}) takes the simpler form
\begin{equation}
\label{bilenbar}
\dfrac{d}{dt}\left( \sum\limits_{i=1}^m {\dot q}_i \dfrac{\partial {\cal L}^*}{\partial {\dot q}_i}-{\cal L}^*\right)=
-\dfrac{\partial {\cal L}^*}{\partial t}+\sum\limits_{\nu=1}^k
\dfrac{\partial T}{\partial {\dot q}_{m+\nu}}\dfrac{\partial \alpha_\nu}{\partial t}.
\end{equation}
\end{prop}

\noindent
{\bf Proof}. Owing to (\ref{baralphaalpha}) the terms with $\dfrac{\partial {\cal L}^*}{\partial q_{m+\nu}}$ in (\ref{bilen}) are null. Moreover, the definition (\ref{baralpha}) shows ${\overline B}_\nu=\dfrac{\partial {\overline \alpha}_\nu}{\partial t}=
\dfrac{\partial \alpha_\nu}{\partial t}$, hence
$\sum\limits_{\nu=1}^k {\overline B}_\nu \dfrac{\partial T}{\partial {\dot q}_{m+\nu}}=
\sum\limits_{\nu=1}^k\dfrac{\partial {\overline \alpha}_\nu}{\partial t} \dfrac{\partial T}{\partial {\dot q}_{m+\nu}}$ and (\ref{bilenbar}) is proved. $\quad\square$

\begin{rem}
The presence of the rheonomic contributions on the right side of equality (\ref{bilenlin}) is easily explainable: $\dfrac{\partial {\cal L}^*}{\partial t}$ can be not null either because of the non--stationarity of the holonomic constraints (hence $\dfrac{\partial T^*}{\partial t}\not =0$)  
or because of the presence of $t$ in the forces (then $\dfrac{\partial U}{\partial t}\not =0$). On the other hand, the possible non--stationarity of the nonholonomic constraints (\ref{constrexpllin}) gives rise to the terms containing $\dfrac{\partial \alpha_\nu}{\partial t}$.
\end{rem}

\subsection{Linear nonholonomic constraints}

\noindent
A significant circumstance of validity of assumption is the case of linear nonholonomic constraints of the form (\ref{constrexpllin}). Indeed, the linear functions $\sum\limits_{j=1}^m \alpha_{\nu, j}(q_1, \dots, q_n, t){\dot q}_j$, $\nu=1, \dots, k$, are homogeneous functions of degree $1$ w.~r.~t.~${\dot q}_1$, $\dots$, ${\dot q}_m$ and Lemma 1 is applicable. By virtue of (\ref{baralphaalpha}) the coefficients (\ref{barbnu}) are ${\overline B}_\nu = \sum\limits_{i=1}^m {\dot q}_i
\dfrac{\partial \alpha_{\nu,i}}{\partial t}$, hence the energy balance (\ref{bilenbar}) takes the form 

\begin{equation}
\label{bilenlin}
\dfrac{d}{dt}\left( \sum\limits_{i=1}^m {\dot q}_i \dfrac{\partial {\cal L}^*}{\partial {\dot q}_i}-{\cal L}^*\right)=
-\dfrac{\partial {\cal L}^*}{\partial t}+\sum\limits_{\nu=1}^k\sum\limits_{i=1}^m {\dot q}_i
\dfrac{\partial T}{\partial {\dot q}_{m+\nu}}\dfrac{\partial \alpha_{\nu,i}}{\partial t}.
\end{equation}
Concerning the linear stationary case, the balance (\ref{bilenlin}) assumes the form pertinent to holonomic systems:

\begin{cor}
For a system such that $\dfrac{\partial {\cal L}^*}{\partial t}=0$ and submitted to stationary linear kinematic constraints (\ref{constrexpllin}) with $\alpha_{\nu,j}=\alpha_{\nu,j}(q_1, \dots , q_n)$ for each $\nu=1, \dots, k$ and $j=1, \dots, m$, the quantity 
$I(q_1, \dots, q_m, {\dot q}_1, \dots, {\dot q}_m)=\sum\limits_{i=1}^m {\dot q}_i \dfrac{\partial {\cal L}^*}{\partial {\dot q}_i}-{\cal L}^*$ is conserved.
Assuming for $T^*$ the form (\ref{tmobile}) with $g_{i,j}$, $b_i$ and $c$ not depending on $t$ explicitly,  $i,j=1, \dots, n$,  
the constant of motion is 
\begin{equation}
\label{intjac}
I=
\sum\limits_{r,s}^m 
\left(\frac{1}{2} g_{r,s}+\frac{1}{2}\sum\limits_{\nu, \mu=1}^k 
 g_{m+\nu, m+\mu} \alpha_{\nu,r} \alpha_{\mu,s}+
\sum\limits_{\nu=1}^k g_{r,m+\nu} \alpha_{\nu, s}
\right){\dot q}_r {\dot q}_s-U-c.
\end{equation}
\end{cor}
The conserved quantity, examined in \cite{neimark} as well as in other textbooks, is the generalized energy integral, or Jacobi integral, of the Lagrangian ${\cal L}^*$. 

\begin{exe}
A very simple model for linear constraints is a point in a pair of points $P_1$ and $P_2$ constrained on vertical plane, keeped at constant distance $2\ell$ one from the other and moving in a way such that the middle point's velocity is along the direction $P_1-P_2$ (\cite{gant}). 
By employing the angle $q_1$ that $P_2-P_1$ forms with the $x$--axis and the coordinates $(q_1, q_2)$ of the middle point, the holonomic setup is $x_1 =q_1+ \ell \cos q_3$, $y_1=q_2+\ell \sin q_3$, $x_2=q_1-\ell \sin q_3$, $y_2=q_2-\ell \sin q_3$. The linear nonholonomic constraint is
${\dot q}_2\sin q_1 - {\dot q}_3\cos q_1=0$ which can be made explicit by writing ${\dot q}_3= {\dot q}_2\tan q_1$. The function (\ref{tmobile}) is $T^*=\frac{1}{2}(M_1+M_2)
\left( \ell^2 {\dot q}_1^2 + (1+\tan^2 q_1) {\dot q}_2^2 \right)$ and, assuming that the $y$--axis points at the upward vertical direction, (\ref{lstar}) is ${\cal L}^*=T^* -(M_1+M_2g)q_3$. The constant of motion (\ref{intjac}) reads 
$I=\frac{1}{2}(M_1+M_2)
\left( \ell^2 {\dot q}_1^2 + (1+\tan^2 q_1) {\dot q}_2^2\right) +(M_1+M_2g)q_3$.
	
\noindent	
We incidentally remark that the model is different form the one corresponding to the requests 
$\overline{P_1P_2}=\ell$, $|{\dot P}_1|=|{\dot P}_2|$ (equidistant points and equal intensity of velocities, see \cite{zek2}): actually, in that case the restrictions yield to the nonlinear condition ${\dot q}_1 ({\dot q}_2\sin q_1 - {\dot q}_3\cos q_1)=0$.
\end{exe}

\begin{exe}
We implement (\ref{intjac}) for a frequently proposed system (\cite{zek3}, \cite{benenti}), consisting in two material points $(P_1, M_1)$ and $(P_2, M_2)$ on a plane whose velocities are both orthogonal to the line joining $P_1$ with $P_2$. The linear kinematic constraints for $P_1\equiv (x_1, y_1)$ and $P_2\equiv (x_2, y_2)$ are 
$$
{\dot x}_1(x_2-x_1)+{\dot y}_1 (y_2-y_1)=0, \quad {\dot x}_2(x_2-x_1)+{\dot y}_2 (y_2-y_1)=0.
$$
Setting $(q_1, q_2, q_3, q_4)=(x_1, x_2, y_1, y_2)$ the explicit kinematic restrictions are 
${\dot q}_3=\frac{q_2-q_1}{q_3-q_4}{\dot q}_1$, ${\dot q}_4=\frac{q_2-q_1}{q_3-q_4} {\dot q}_2$ and 
one has, with respect to (\ref{gxi}), (\ref{constrexpllin}) (in this case it is $m=k=2$):
$$
g_{1,1}=g_{3,3}=m_1, \quad g_{2,2}=g_{4,4}=m_2, \quad g_{i,j}=0\; \;\textrm{for}\;i\not =j, \quad
\alpha_{1,1}=\alpha_{2,2}=\frac{q_2-q_1}{q_3-q_4}, \;\;\alpha_{2,1}=\alpha_{2,2}=0
$$
so that (\ref{intjac}) gives the conserved quantity
$$
I(q_1, q_2, q_3, q_4,{\dot q}_1, {\dot q}_2)=
\frac{1}{2}(m_1{\dot q}_1^2 +m_2{\dot q}_2^2)\left(1+\left(\frac{q_2-q_1}{q_3-q_4}\right)^2\right)
-U(q_1, q_2, q_3, q_4)
$$
where $U$ is calculated (according to the dynamics of the model) by means of (\ref{potu}).

\noindent
We judge wortwhile to remark that, imaging a transfer from Example 5 to Example 10 (both regard parallel velocities) made by specifiying a request on the velocity vectors, entails a very different structure of the nonholonomic constraints (nonlinear towards linear); such an argument is 
properly highlighted in (\ref{benenti}).

\end{exe}

\begin{exe}
A further model, recurrent in literature, consists in considering two material points $P_1$ and $P_2$ on a plane whose velocities are perpendicular; furthermore, the velocity of $P_1$ is orthogonal to the straight line joining $P_1$ with $P_2$. The model is different from the one discussed in Example 2, where only the nonlinear nonholonomic constraint (\ref{velperp}) is considered: in the present case the linear kinematic condition 
${\cal C}_1:$ ${\dot x}_1 (x_2-x_1)+{\dot y}_1 (y_2 - y_1)=0$ is added, according to the imposed restriction. The same restriction, holding (\ref{velperp}), can be formulated also stating that the velocity of $P_2$ is parallel to the line joining the two points, that is ${\cal C}_2:$ ${\dot x}_2 (y_2-y_1)-{\dot y}_2 (x_2-x_1)=0$. 

\noindent
We can easily show that either (\ref{velperp}) coupled with ${\cal C}_1$ or (\ref{velperp}) coupled with ${\cal C}_2$ are equivalent to the pair of linear constraints ${\cal C}_1$, ${\cal C}_2$. Indeed, it suffices to write the three conditions as $\frac{{\dot x}_1}{{\dot y}_1}=-\frac{{\dot y}_2}{{\dot x}_2}$,  
$\frac{{\dot x}_1}{{\dot y}_1}=-\frac{y_2-y_1}{x_2-x_1}$, $\frac{{\dot y}_2}{{\dot x}_2}=\frac{y_2-y_1}{x_2-x_1}$ to conclude. 
Thus, examining the linear nonholonomic system ${\cal C}_1$, ${\cal C}_2$ we set $(q_1, q_2, q_3, q_4)=(x_1, x_2, y_1, y_2)$ in order to write (\ref{constrexpllin}) as 
${\dot q}_3=-\frac{q_2-q_1}{q_4-q_3}{\dot q}_1$, ${\dot q}_4=\frac{q_4-q_3}{q_2-q_1} {\dot q}_2$ 
and, in a very similar way to the previous example, one finds the conservation of the quantity (\ref{intjac})
$$
I(q_1, q_2, q_3, q_4,{\dot q}_1, {\dot q}_2)=
\frac{1}{2}M_1{\dot q}_1^2\left( 1+\left(\frac{q_2-q_1}{q_4-q_3}\right)^2\right)+
\frac{1}{2}M_2{\dot q}_2^2\left( 1+\left(\frac{q_4-q_3}{q_2-q_1}\right)^2\right)
-U(q_1, q_2, q_3, q_4).
$$
\end{exe}

\subsection{Nonlinear constraints: homogeneous quadratic functions}

\noindent
A frequently encountered subcategory of nonholonomic constraints (\ref{constr}) encompasses restrictions of the form
\begin{equation}
\label{homquadr}
\sum\limits_{i,j=1}^n \gamma_{i,j}^\nu(q_1, \dots, q_n, t) {\dot q}_i {\dot q}_j =0, \quad \nu=1, \dots, k.
\end{equation}
Assuming that the explicit form which singles out ${\dot q}_{m+1}$, $\dots$, ${\dot q}_n$ can be achieved, 
(\ref{constrexpl}) takes the form
\begin{equation}
\label{alphaquadr}
{\dot q}_{m+\nu}=\dfrac{\sum\limits_{r,s=1}^m \gamma_{r,s}^\nu (q_1, \dots, q_n, t){\dot q}_r {\dot q}_s}{\sum\limits_{i=1}^m \beta_i^\nu (q_1, \dots, q_n,t){\dot q}_i}=\alpha_\nu (q_1, \dots, q_n, {\dot q}_1, \dots, {\dot q}_m, t), \qquad \nu=1, \dots, k
\end{equation}
for appropriate functions $\beta_i^\nu$.
The functions $\alpha_\nu$ of (\ref{alphaquadr}) are homogeneous of degree 1 w.~r.~t.~${\dot q}_1$, $\dots$, ${\dot q}_m$, hence the energy balance which pertains to such systems is (\ref{bilenbar}), where the functions $\gamma_{r,s}^\nu$ and $\beta_i^\nu$ will appear; in the stationary case the energy is conserved. In a natural way, constraints (\ref{homquadr}) appear when the restrictions concern parallelism or orthogonality of the velocities, or the assignment of equal intensity of the velocities. 
A typical and simple instance takes into consideration two points $P_1$, $P_2$, for which the three restrictions ${\dot P}_1\wedge {\dot P}_2=0$, ${\dot P}_1 \cdot {\dot P}_2=0$, $|{\dot P}_1|^2 = 
|{\dot P}_2|^2$ respectively read, in  cartesian coordinates, 
$$
\begin{array}{ll}
{\dot x}_1{\dot y}_2 - {\dot y}_1 {\dot x}_2, \;\;{\dot x}_1 {\dot z}_2 - {\dot z}_1 {\dot x}_2=0 & \textrm{parallelism}\\
{\dot x}_1 {\dot x}_2+{\dot y}_1{\dot y}_2+{\dot z}_1{\dot z}_2=0 & \textrm{orthogonality}\\
{\dot x}_1^2 +{\dot y}_1^2 +{\dot z}_1^2 - {\dot x}_2^2 -{\dot y}_2^2 -{\dot z}_2^2=0 & \textrm{same norm of velocity}
\end{array}
$$
belonging to class (\ref{homquadr}). 

\noindent
The already performed Examples 1,4 (nonholonomic pendulum), 2 (perpendicular velocities),
5 (parallel velocities), 7 (equal velocity in norm) are part of the type (\ref{homquadr}) of restrictions.
It is worthwhile to stress, as some of the Examples revealed, that adding further resctrictions or 
specifying mechanisms regarding the nonholonomic constraints (\ref{homquadr}) may modify totally the typology of the restrictions (see Example 9): in this sense, the equivalence in realizing phisically nonholonomic restrictions by means of either linear or nonlinear equations claimed in (\ref{zek}) needs to be debated.

\subsection{Linear affine constraints}

\noindent
Finally, a special situation concerns the nonholonomic systems with affine constraints, which can be assumed of the form (not encompassed by (\ref{constrexpllin}))
\begin{equation}
\label{constraff}
\alpha_\nu=\sum\limits_{j=1}^m a_{\nu,j}(q_1,\dots, q_n, t){\dot q}_j+c_\nu(t), \quad \nu=1,\dots, k
\end{equation}
with $c_\nu$ non zero function. 
In that case ${\bar \alpha}_\nu = \alpha_\nu - c_\nu$, so that  
${\overline B}_\nu=\dfrac{\partial {\overline \alpha}_\nu}{\partial t}$ and (\ref{bilen}) is
\begin{equation}
\label{bilenaff}
\dfrac{d}{dt}\left( \sum\limits_{i=1}^m {\dot q}_i \dfrac{\partial {\cal L}^*}{\partial {\dot q}_i}-{\cal L}^*\right)
+c_\nu\dfrac{\partial {\cal L}^*}{\partial q_{m+\nu}}-\sum\limits_{\nu=1}^k
\left( \sum\limits_{j=1}^m \dfrac{\partial \alpha_{\nu,j}}{\partial t}{\dot q}_j+{\dot c}_\nu\right) \dfrac{\partial T}{\partial {\dot q}_{m+\nu}}
=
-\dfrac{\partial {\cal L}^*}{\partial t}.
\end{equation}
The stationary case $\alpha_{\nu,j}(_1, \dots, q_n, {\dot q}_1, \dots, {\dot q}_m)$, $c_\nu$ constant 
provides the conservation of the quantity in round brackets (energy), whenever ${\cal L}^*$ does not depend explicitly on $t$ and the forces verify special pro\-per\-ties, as it is described in the following

\begin{exe}
For a particle in ${\Bbb R}^3$ with mass $M$ and submitted to the affine constraint $ax {\dot y}+b {\dot x}y +c -{\dot z}=0$, $c\not=0$, we set $(q_1, q_2, q_3)=(x,y,z)$ so that (\ref{constraff}) writes
${\dot q}_3=a q_1{\dot q}_2+b q_2{\dot q}_1  +c$ (in this case $m=2$ and $k=1$); moreover 
$${\cal L}^*=\frac{1}{2}M [(1+a^2 q_1^2){\dot q}_1^2+(1+b^2q_2^2){\dot q}_2^2+2ab q_1q_2 {\dot q}_1{\dot q}_2+2c (aq_1{\dot q}_2 + b q_2 {\dot q}_1)]-V(q_1, q_2, q_3).
$$
Whenever $V=V(q_1, q_2)$, (\ref{bilen}) supplies the conservation of the quantity
$$
\frac{1}{2}M({\dot q}_1^2+{\dot q}_2^2)-\frac{1}{2}M (aq_1 {\dot q}_2 +b q_2 {\dot q}_1 +c)(aq_1 {\dot q}_2 +b q_2 {\dot q}_1 -c).
$$
\end{exe}

\noindent
A necessary and sufficient condition (in terms of geometrical properties of the constraint manifold) in order that (\ref{bilen}) provides the energy integral is discussed and proved in \cite{fas}.

\section{Conclusions}

\noindent
The equations of motions for nonholonomic nonlinear systems are presented in the double version 
(\ref{vnl2}) and (\ref{vnl}), each of them showing advantageous points. The context of nonlinear restrictions leads us to make use of the generalized velocities as kinematic variables (instead of quasi--coordinates), thus favouring the extension of the Voronec's method for linear kinematic constraints to the nonlinear case.

\noindent
The calculation of the power of the forces, by way of the equations of motion, generates the formula (\ref{bilen}), showing in an unified and consistent way the rate of change in time of the energy espressed by the independent velocities (the function in round brackets in (\ref{bilen})) in terms of the contributions due to the constraints forces (by means of the coefficients ${\overline \alpha}_\nu$ and ${\overline B}_\nu$) and of the possible explicit dependence of the restriction or of the forces on time (term on the right hand side).

\noindent
The arrangement of the energy balance equation is suitable in order to identify the type of nonholonomic constraints exhibiting ${\overline \alpha}_\nu=\alpha_\nu$, which is the key point for the purposes of inferring the first integral of energy.
At the same time, in some special modellistic circumstances the setup (\ref{bilen}) shows directly the appropriate simplifcations, as in the case (\ref{bilenchap}) or in the linear case.
Several examples of simple but meaningful systems have been performed. 

\noindent
The condition (\ref{baralphaalpha}) fits for kinematic constraints whose explicit formulation (\ref{constrexpl}) is a homogeneous function of degree $1$: the models frequently adopted in literature and accessible for implementations concerning special restrictions on the velocities (parallelism, equal norms, orthogonality)  fall in this category. Theoretically, it is easy to extend the category to restrictions (\ref{constr}) which are homogeneous functions of arbitrary degree w.~r.~t.~the generalized velocities ${\dot q}_1$, $\dots$, ${\dot q}_n$; however, from the experimental point of view this may produce a not concrete realisation.

\noindent
Beyond the mere aspect of the balance of energy, the subject nonlinear constraints presents interesting questions somehow unexplored in literature and sometimes misleading. These aspects, just mentioned in the paper, concern the equivalence of linear kinematics models with nonlinear restrictions, the correctness of mergering part of the constraints giving rise to new conditions (tipically, two linear constraints are joined to form a quadratic condition); the Hamel--Appell example of a system with nonlinear nonholonomic constraint obtained by linear kinematic condition is a point of reference in this sense.

\noindent
A sistematic procedure for readily comparing the equations of motion whenever different sets of independent velocities are selected is also not secondary in order to take into the right consideration the local use of (\ref{constrexpl}). 

\noindent
The just mentioned open points are the next purpose of the research on nonlinear kinematic constraints.

\end{document}